\input amssym.def
\input amssym
\magnification=1200
\parindent0pt
\hsize=16 true cm
\baselineskip=13  pt plus .2pt
$ $

\def\K{(\Gamma',{\cal G'})}
\def\G{(\Gamma,{\cal G})}
\def\Z{{\Bbb Z}}
\def\D{{\Bbb D}}

\def\H{{\cal H}}

\centerline {\bf On finite groups acting on a connected sum of 3-manifolds $S^2
\times S^1$}

\medskip

\centerline {by}

\medskip

\centerline {{\bf Bruno P. Zimmermann}  (Trieste)}

\bigskip

{\bf Abstract.} Let $H_g$ denote the closed 3-manifold obtained as the
connected sum of $g$ copies of $S^2 \times S^1$, with free fundamental group of
rank $g$.  We prove that, for a finite group $G$ acting on $H_g$ which induces
a faithful action on the fundamental group, there is an upper bound for the
order of $G$ which is quadratic in $g$, but that there does not exist a linear
bound in $g$. This implies then a Jordan-type bound for arbitrary finite group
actions on $H_g$ which is quadratic in $g$.  For the proofs we develop a
calculus for finite group-actions on $H_g$, by codifying such actions by
handle-orbifolds and finite graphs of finite groups.

\bigskip


{\sl Key words and phrases:}  3-manifold, connected sum of 3-manifolds $S^2
\times S^1$, bound on finite group action,  Jordan-type bound

\bigskip

{\bf 1. Introduction. } All finite group-actions in the present paper will be faithful,
smooth and orientation-preserving, and all manifolds and orbifolds will be
orientable. For the case of surfaces, the famous Hurwitz bound states that the
order of a finite group acting on a closed surface of genus $g > 1$ is bounded
above by $84(g-1)$. In a similar spirit, the order of a finite group acting on
a 3-dimensional handlebody $V_g$ of genus $g>1$ is bounded by $12(g-1)$ ([Z1],
[MMZ, Theorem 7.2]).

\medskip

We consider finite groups $G$ acting on the connected sum
$H_g = \sharp_g (S^2 \times S^1)$ of $g$ copies of $S^2 \times S^1$; we will call
$H_g$ a  {\it closed handle} or just a {\it handle} of genus $g$ in the following.
Similar as for a handlebody of genus $g$, the fundamental group of a closed handle $H_g$ is
the free group $F_g$ of rank $g$.  For a handlebody $V_g$ of genus
$g>1$, a finite group
$G$ acting on $V_g$ acts faithfully on its fundamental group, meaning that it induces an
injection  $G \to {\rm Out} \, \pi_1(V_g)$ into the outer automorphism group ${\rm Out} \,
\pi_1(V_g)$ (see Proposition 2), and hence induces a subgroup of order at most $12(g-1)$ of
${\rm Out} \, F_g$.  The maximum order of a general finite subgroup of ${\rm Out} \,
F_g$ instead is  $2^g g!$ for $g > 2$ and 12 for
$g=2$ ([WZ])  (based on the result in [Z2]  that each finite subgroup of
${\rm Out} \, F_g$ can be induced by an action of the group on a finite graph, and then
also on some handlebody of sufficiently high dimension).

\medskip

Since $H_g$  admits $S^1$-actions (see [R]), it admits finite cyclic group-actions of
arbitrarily large order acting trivially on the fundamental group (i.e., inducing the
trivial homomorphism to ${\rm Out} \, F_g$), so in contrast to the situation for
handlebodies of genus
$g>1$ there is no upper bound for the orders of finite groups
$G$ acting on a closed handle $H_g$.  However, let $G_0$ denote the normal subgroup of all
elements of $G$ acting trivially on $\pi_1(H_g)$; then $G_0$
is cyclic for $g>1$ (Proposition 2), the quotient $H_g/G_0$ is again a closed
handle of the same genus $g$, and the factor group
$G/G_0$ acts faithfully on the fundamental group of the quotient $H_g/G_0 \cong H_g$. Hence
one is led to consider actions of finite groups $G$ on
$H_g$ which act faithfully on the fundamental group; in section 4, we prove the following:

\bigskip

{\bf Theorem 1.}  {\sl  Let $G$ be a finite group acting on a closed
handle $H_g$ of genus $g$ such that the induced action on the fundamental group is
faithful. Then, for $g \ge 15$, there is the quadratic bound
$$|G| \; \le  \; 24g(g-1).$$
If $G$ is cyclic then, for $g \ge 3$,
$$|G| \; \le  \; 4(g-1)^2.$$
There does not exist a linear bound in $g$ for the order of $G$.}

\bigskip

By the above, Theorem 1 implies the following Jordan-type bound for arbitrary finite groups
acting on a closed handle (i.e., not necessarily faithful on the fundamental group):

\bigskip

{\bf Corollary.}  {\sl  Let $G$ be a finite group acting on a closed handle
$H_g$ of genus $g > 1$. Then $G$ has a cyclic normal subgroup $C$ (the subgroup
acting trivially on the fundamental group) such that, for $g \ge 15$, the order
of $G/C$ is bounded above by  $24g(g-1)$.}

\medskip

We note that by the classical Jordan bound each finite subgroup $G$ of a linear
group GL$(n,\Bbb C)$ has a normal abelian subgroup $A$ such that the order of
$G/A$ is bounded  by a constant depending only on $n$; by [C], for $n \ge 71$
the optimal bound here is $(n+1)!$, realized by the symmetric group $S_{n+1}$
occurring as a subgroup of GL$(n,\Bbb C)$.

\medskip

After the various linear bounds for finite group-actions on surfaces and
handlebodies, Theorem 1 seems to present the first instance of a quadratic bound in such
a situation. There remains the problem to determine the optimal quadratic
bound (the optimal coefficient for $g^2$), for both cyclic and arbitrary finite groups; we
will construct some explicit examples of finite cyclic and finite group-actions on
closed handles in section 3 which seem to come close to these optimal bounds.

\medskip

In order to prove Theorem 1 we shall develop in section 2 a calculus for
finite group-actions on closed handles $H_g$ (see Theorem 2), in analogy with the theory
of finite group-actions on handlebodies $V_g$ in [MMZ] (see also [Z3], and [MZ] for
applications to group-actions of large order on handlebodies and closed 3-manifolds). This
uses the language of handle-orbifolds and of finite graphs of finite groups which codify the
quotient orbifolds  $H_g/G$.

\medskip

We note that the maximum order of a finite cyclic group acting on a
closed surface of genus  $g>1$ is
$4g+2$, for a handlebody of genus $g>1$ it is
$2g+2$ if $g$ is even, and $2g-2$ if $g$ is odd ([MMZ]). For general finite cyclic
subgroups of ${\rm Out} \, F_g$ instead, it is proved in [LN] and [B] that the maximum order
behaves   approximatively like the Landau estimate
${\rm exp} (\sqrt {g \; {\rm log} \, g})$ for  the maximum orders of elements of
the symmetric group of degree $g$. The maximum orders of the finite cyclic
subgroups of Aut($F_g$) and  Out($F_g$) can be found in [LN, Table 1], for $g \le 300$.

\medskip

We  close the introduction with the question of what happens for the
higher-dimensional analogues  $\sharp_g (S^d \times S^1)$ of a closed handle
$H_g$, for $d > 2$ (all with free fundamental group of rank $g$). In analogy
with Theorem 1, is there a polynomial bound in $g$ for the order of $G$, and of
what degree? (See also [Z2] for a discussion of finite group-actions on
higher-dimensional analogues of handlebodies.)

\bigskip

{\bf 2. Handle orbifolds and associated finite graphs of finite groups.}
Let $G$ be a finite group acting faithfully and orientation-preservingly on a
handle
$H_g = \sharp_g (S^2 \times S^1)$ of genus $g$.  Denoting by $E$ the group generated by
all lifts of elements of $G$ to the universal covering of
$H_g$ and by $F_g$ the normal subgroup of covering transformations, we have a
group extension  $ 1 \to F_g  \to  E \to G \to 1$ which belongs to the abstract
kernel $\;G \to  {\rm Out} \, \pi_1(H_g) \cong {\rm Out} \, F_g$ induced by the
action of $G$ on $\pi_1(H_g)$.  Using the equivariant sphere
theorem (see [MSY] for an approach by minimal surface techniques, and [Du] and [JR]
for purely topological-combinatorial proofs), we will associate to the action of
$G$ a  handle-orbifold
$\H$ and a finite graph of finite groups $\G$ whose fundamental group
$\pi_1\G$ is isomorphic to the extension $E$.

\medskip

By the equivariant sphere theorem there exists an embedded, homotopically nontrivial
2-sphere  $S^2$ in $H = H_g$ such that $x(S^2) = S^2$ or
$x(S^2) \cap S^2 = \emptyset$ for all $x \in G$. We cut $H$ along the system
of disjoint 2-spheres $G(S^2)$, by removing the interiors of $G$-equivariant regular
neighbourhoods $S^2 \times [-1,1]$ of each of these 2-spheres, and call each
of these regular neighbourhoods $S^2 \times [-1,1]$ a 1-handle. The
result is a collection of 3-manifolds with 2-sphere boundaries, with an induced
action of $G$. We close each of the 2-sphere boundaries by a 3-ball and extend
the action of $G$ by taking the cone over the  center of each of
these 3-balls, so $G$ permutes these 3-balls and their centers. The result is a
finite collection of closed handles of lower genus on which $G$ acts.
Applying inductively the procedure of cutting along 2-spheres, we finally
end up with a finite collection of 3-spheres which we call
0-handles (these are just the closed handles of genus 0). Note that the construction
gives a finite graph $\tilde \Gamma$ on which $G$ acts whose vertices correspond to the
0-handles and whose edges to the 1-handles.

\medskip

On each 3-sphere (0-handle) there are finitely many points which are the
centers of the attached 3-balls (their boundaries are the 2-spheres along which
the 1-handles are attached).  For each of these 3-spheres we consider its
stabilizer $G_v$ in $G$. By the recent geometrization of finite group-actions
on 3-manifolds following Thurston and Perelman, we can assume that the action
of a stabilizer $G_v$  on the corresponding 3-sphere is standard, i.e.
orthogonal; we call such a quotient $S^3/G_v$ a 0-handle orbifold (we
note that the geometrization is not really essential for the construction and
its applications, it just says that each 0-handle orbifold is standard).
Similarly, considering the stabilizers $G_e$ in $G$ of the 1-handles $S^2
\times [-1,1]$, we can assume that each stabilizer $G_e$ preserves the product
structure of $S^2 \times [-1,1]$. If some element of a stabilizer $G_e$ acts as
a reflection on [-1,1], we split the 1-handle into two 1-handles by introducing
a new 0-handle obtained form a small regular neighbourhood $S^2 \times
[-\epsilon,\epsilon]$ of $S_2 \times \{0\}$ by closing up with two 3-balls.
Hence we can assume that each stabilizer $G_e$ of a 1-handle $S^2 \times
[-1,1]$ does not interchange its two boundary 2-spheres; equivalently, $G$ acts
without inversions on the graph $\tilde \Gamma$. We call such a quotient
$(S^2 \times [-1,1])/G_e \cong (S^2/G_e) \times [-1,1]$ a 1-handle
orbifold.

\medskip

As result, the quotient orbifold $\H = H/G$ is obtained from a finite collection of
0-handle orbifolds $S^3/G_v$ by removing  the interiors of disjoint 3-ball
neighbourhoods of finitely many points and attaching 1-handle orbifolds along the
resulting 2-sphere boundaries (respecting singular sets and branching orders as well as
orientations); we call such a structure a {\it closed handle-orbifold} or just a {\it
handle-orbifold}. Summarizing, we have:

\bigskip

{\bf Proposition 1.}  {\sl The quotients of closed handles $H_g$ by finite group-actions
have the structure of closed handle-orbifolds.}

\bigskip

We note also the following easy but crucial

\medskip

{\bf Observation.}  If around a point of a  0-handle $S^3$ of $H$ a 1-handle is attached
then the stabilizer $G_e$ of the 1-handle is exactly the subgroup of the
stabilizer $G_v$ of the 0-handle which fixes the considered point (since otherwise some
larger subgroup of $G_v$ would stabilize the 1-handle).

\bigskip

To each handle-orbifold $\H = H/G$ is associated a graph of groups $\G$ in a
natural way. The underlying graph  $\Gamma$ is just the quotient $\tilde
\Gamma/G$. The vertices (resp. edges) of the graph $\Gamma$ correspond to
0-handle orbifolds $S^3/G_v$ (resp. 1-handle orbifolds $(S^2/G_e) \times
[-1,1]$), and to each vertex (resp. edge) is associated the corresponding
stabilizer $G_v$ (resp. $G_e$) (choosing an isomorphic lift of a maximal tree
of $\Gamma$ to $\tilde \Gamma$, and then also lifts of the remaining edges).
In particular, the vertex groups $G_v$ of $\G$ are isomorphic to finite
subgroups of the orthogonal group SO(4), and the edge groups $G_e$ to finite
subgroups of SO(3). We can also assume that the graph of groups has no trivial
edges, i.e. no edges with two different vertices such that the edge group
coincides with one of the two vertex groups (by collapsing the edge, i.e.
amalgamating the two 0-handles into a single 0-handle). We say that such a
handle-orbifold $\H$ and the associated graph of groups is in {\it normal
form}.

\medskip

We have associated to each handle-orbifold $\H = H/G$ a graph of
groups $\G$ in normal form. By the orbifold version of Van Kampen's theorem (see [HD]),
the orbifold fundamental group $\pi_1(\H)$ is isomorphic to the  fundamental
group $\pi_1\G$ of the graph of groups $\G$ (which is the iterated free product
with amalgamation and HNN-extension of the vertex groups over the edge groups,
starting with a maximal tree in $\Gamma$; see [Se], [ScW] or [Z4] for the standard
theory of graphs of groups, their fundamental groups  and the connection with groups
acting on trees and graphs).  We have also a canonical surjection $\phi: \pi_1(\H) \cong
\pi_1\G \to G$, injective on vertex and edge groups, whose kernel is
isomorphic to the fundamental group
$\pi_1(H)
\cong F_g$ of the handle
$H$, and the group extension
$$ 1 \to \pi_1(H) \cong F_g  \to  \pi_1({\H})  \cong \pi_1\G  \to G \to 1$$
is equivalent to the group extension $1 \to F_g \to E  \to G$. In particular, $H$
is the orbifold covering of $\H$ associated to the kernel of the surjection $\phi$.

\bigskip

Conversely, suppose we have a finite graph of finite groups $\G$  associated to a
handle-orbifold $\H$ and a surjection
$\phi: \pi_1\G  \to G$ onto a finite group $G$ which is injective on the vertex groups.
Then the orbifold covering of $\H$ associated to
the kernel of $\phi$ is a closed handle $H_g$ of some genus $g$ on which $G$ acts as
the group of covering transformations. The genus $g$ can be computed as
follows. Denoting by
$$\chi\G = \sum {1 \over |G_v|} - \sum {1 \over |G_e|}$$
the  Euler characteristic of the graph of groups $\G$ (the sum is taken over all
vertex groups $G_v$ resp. edge groups $G_e$ of $\G$), we have
$$g-1 =  -\chi\G \; |G|$$
(see [ScW], [Z4]).

\medskip

Finally, we note that the induced action of $G$ on the fundamental group of $H$ is
effective (faithful) if and only if the corresponding group extension
$1 \to F_g  \to E  \to G \to 1$ is effective (i.e., by considering $F_g$ as a subgroup of
$E$, the homomorphism
$G \to {\rm Out} \, F_g$ induced by conjugation of $F_g$ by preimages in $E$ of
elements in $G$ is injective). It is easy to see that this is the case if and only if the
extension group $E \cong \pi_1\G$ has no nontrivial finite normal subgroups: If a
preimage $e \in E$ of an element $x \in G$ induces by conjugation an inner automorphism
of $F_g$, then another preimage of $x$ in $E$ induces the trivial or identity automorphism
of $F_g$; since a power of $e$ lies in $F_g$ and the center of $F_g$ is trivial,  $e$ must
have finite order, and clearly the subgroup of elements of $E$ acting by conjugation
trivially on $F_g$ is a finite normal subgroup of $E$).

\medskip

Summarizing, we have:

\bigskip

{\bf Theorem 2.}  {\sl  A finite group $G$ acts on a closed handle $H_g$ of genus $g$ if and
only if there is a finite graph of finite groups $\G$ in normal form associated to a
handle-orbifold $\H$, and a surjection $\phi: \pi_1\G \to G$ which is
injective on the vertex groups such that
$$g =  -\chi\G \; |G| + 1.$$
The induced action of $G$ on the fundamental group of $H_g$ is faithful if and
only if $\pi_1\G$ has no nontrivial finite normal subgroups.}

\medskip

See section 3 for some significant examples.  As noted in the introduction, also the
following holds.

\bigskip

{\bf Proposition 2.}  {\sl i)  Let $G$ be a finite group acting on closed handle $H_g$ of
genus $g>1$. Then the normal subgroup $G_0$ of all elements of $G$  inducing a
trivial action on the fundamental group  is  cyclic, and the quotient $H_g/G_0$ is again
homeomorphic to a closed handle of genus $g$.

\smallskip

ii)  Let $G$ be a finite group acting faithfully on a handlebody $V_g$. If $g>1$ then the
induced action of $G$ on the fundamental group is faithful.}

\bigskip

{\it Proof.} i)  Consider the action of $G$ on a graph $\tilde \Gamma$
associated to a handle-decomposition of $H_g$ as before; in particular, $\tilde
\Gamma$ has no vertices of degree 1. Since $g>1$, by [Z5, Lemma 1] (considering
homology and the Hopf trace formula) the action of $G_0$ on $\tilde \Gamma$ has
to be trivial (or by a direct combinatorial argument). Hence $G_0$ maps each
1-handle $S^2 \times [-1,1]$ and each 0-handle $S^3$ to itself and is
isomorphic to a subgroup of SO(3). Moreover $G_0$ has to be a cyclic group
since a non-cyclic group $G_ 0$  would have only two global fixed points in
each 0-handle $S^3$ around which a 1-handle can be attached, so the graph
$\tilde \Gamma$ would be a circular graph and $g = 1$ (a cyclic group instead
has a circle of fixed points in each 0-handle along which arbitrarily many
1-handles can be attached).

\medskip

ii)  This follows similar similar as in the proof of i) from the analogous
theory of finite group-actions on handlebodies, replacing the equivariant
sphere theorem by the equivariant Dehn lemma/loop theorem. Denoting by $B^n$
the closed $n$-ball, the stabilizers of the 1-handles $B^2 \times [-1,1]$  are
finite subgroups of SO(2) now, the stabilizers of the 0-handles $B^3$ are
finite subgroups of SO(3), and each 1-handle is attached along its two boundary
components to the boundary of  one or two 0-handles $B^3$;  since $g>1$, $G_0$
has to be trivial now (alternatively, one  may apply [Z1, Korollar 1.3]).

\bigskip

{\bf 3. Examples.}  We construct first an infinite series of finite cyclic group-actions on
closed handles, faithful on the fundamental group. For large $g$, this realizes the maximum
order for cyclic group-actions which we know at present.

\medskip

For an odd positive integer $a$ which is
divisible by 3, consider $G \cong \Z_n$ where $n = a(a+1)$. Let  $\G$ be the graph of groups
which consists of two edges with edge groups $\Z_a$ and $\Z_{a+1}$, and three
vertices with vertex groups $\Z_{2a}$, $\Z_{3(a+1)}$ (both of valence 1) and
$\Z_{a(a+1)}$ (the middle vertex of valence 2); its fundamental group is the free product
with amalgamation
$$\pi_1\G \;\; \cong  \;\;  \Z_{2a}
*_{\Z_a}\Z_{a(a+1)}*_{\Z_{a+1}}\Z_{3(a+1)};$$
note that $\pi_1\G$ has no nontrivial finite normal subgroups (cf. Theorem 2).

\medskip

Choose an orthogonal action of
$\Z_{a(a+1)}$ on $S^3$ such the subgroups $\Z_a$ and $\Z_{a+1}$ have two
disjoint circles of fixed points, and associate a handle-orbifold $\H$ to
$\G$, with 0-handles  $S^3/\Z_{a(a+1)}$, $S^3/\Z_{2a}$ and $S^3/\Z_{3(a+1)}$ such that
$\pi_1\H \; \cong \pi_1\G$ (cf. the Observation after Proposition 1).

\medskip

There is an obvious surjection $\phi: \pi_1\G \to \Z_n$,
injective on vertex groups, which defines an action of $\Z_n$ on a closed handle $H_g$ (the
orbifold covering of the handle-orbifold $\H$ corresponding to the kernel of $\phi$). Now

$$-\chi \;  = \;  -\chi\G \; = \; {1 \over a} + {1 \over a+1} - {1 \over 2a} - {1 \over
a(a+1)} - {1
\over 3(a+1)} \;  =  \; {7a-3 \over 6a(a+1)}, $$
$$g-1 \; = \;  -\chi \;  n  \; = \;  {7a-3 \over 6}, \hskip 7mm  a \; = \; {6g-3 \over
7}$$ and finally
$$n \; = \; a(a+1) \; = \; {(6g-3) \over 7} \;\; {(6g+4) \over 7}$$
 which is quadratic in $g$. In particular, there cannot exist a linear bound in $g$ for
the order of $G$, proving the final assertion of Theorem 1.

\bigskip \bigskip

Next we discuss an infinite series of actions of non-cyclic groups
on closed handles, faithful on the fundamental group; again this realizes the maximum
order for arbitrary finite group-actions which we know at present, for large values of $g$.

\medskip

For arbitrary finite groups $G$, the vertex groups of a graph of groups $\G$ as
in Theorem 2  are finite subgroups of the orthogonal group SO(4), and the edge
groups are finite subgroups of SO(3). The  orthogonal group SO(4) is isomorphic
to the central product $S^3 \times_{\Z_2} S^3$ of two copies of the unit
quarternions $S^3$ (i.e., with identified centers $\Z_2$); the orthogonal
action on $S^3$ is given by $x \to q_1^{-1}xq_2$, for a fixed pair of unit
quarternions $(q_1,q_2) \in S^3  \times S^3$. There is a 2-fold covering $S^3
\to {\rm SO} (3)$ whose kernel is the central subgroup $\Z_2 = \{\pm 1\}$ of
$S^3$, and the finite subgroups of $S^3$ are exactly the binary polyhedral
groups which are the preimages of the polyhedral groups in SO(3); we denote by
$\D_{4a}^*$ the binary dihedral group of order $4a$ which is the the preimage
of the dihedral group $\D_{2a}$ of order $2a$, with $\D_{4a}^*/\Z_2 \cong
\D_{2a}$.

\medskip

For an integer $a \ge 2$, let $\G$ be the graph of groups consisting of a single edge and
two vertices, with vertex groups
$\; \D_{4a}^* \times_{\Z_2} \D_{4a}^*  \,  \subset  \, S^3 \times_{\Z_2} S^3$ (of order
$8a^2$) and $\, \D_{2a} \times \Z_2$ (of order $4a$) where $\, \D_{2a}$ denotes the
diagonal subgroup of $\, \D_{4a}^* \times_{\Z_2} \D_{4a}^*$ (which is the maximal
subgroup fixing the points $\pm 1 \in S^3$) and $\Z_2$ the central subgroup of $S^3$; the
edge group is the dihedral group $\D_{2a}$. Its fundamental group is the free product with
amalgamation
$$\pi_1\G \; \cong \; (\D_{4a}^* \times_{\Z_2} \D_{4a}^*) \, *_{\D_{2a}}
(\D_{2a} \times {\Z_2}),$$
and there is an obvious surjection, injective on vertex groups,
$$\phi: \pi_1\G \; \to \; \D_{4a}^* \times_{\Z_2} \D_{4a}^*.$$
The order of any finite group $G$ onto which $\pi_1\G$ surjects, injective on
vertex groups, has some order $8xa^2$. Now
$$ -\chi \; = \;  -\chi\G  \; =  \; {1 \over 2a} - {1 \over 8a^2} -
{1 \over 4a} \; = \; {2a-1 \over 8a^2},$$
$$g-1 \; =  \; -\chi \; 8xa^2  \;  = \; x(2a-1), \hskip 5mm g =   x(2a-1) + 1,$$
$${|G| \over g^2} \; =  \; {8xa^2 \over x^2(2a-1)^2 + 1 + 2x(2a-1)} \;  \le \;
{8a^2  \over (2a-1)^2 + 1 + 2(2a-1)} \; = \; {8a^2 \over 4a^2} \; = \; 2.$$

\bigskip

It follows that  $|G| \le 2g^2$, and
$$|G| = 2g^2$$

if and only if $x=1$ and $G \cong
\D_{4a}^* \times_{\Z_2} \D_{4a}^*$.

\bigskip

{\bf 4. Proof of Theorem 1.}   Suppose that the finite group $G$ of order $n$ acts on a
closed handle $H = H_g$ of genus $g>1$, faithfully on the fundamental group. By Theorem 2
there is a finite graph of finite groups
$\G$ in normal form associated to a  handle-orbifold ${\H}$ whose fundamental group
$\pi_1(\H)
\cong
\pi_1\G$ has no nontrivial finite normal subgroups, and a surjection  $\phi: \pi_1\G \to
G$, injective on vertex groups, such that the action of
$G$ on $H$ is given by the orbifold covering of ${\H}$ associated to
the kernel of $\phi$.  The edge groups of $\G$ are finite subgroups of SO(3), that is
cyclic, dihedral, tetrahedral of order 12, octahedral of order 24 or dodecahedral of order
60.   Let $\chi = \chi\G$ denote the Euler characteristic of $\G$; note
that $-\chi > 0$ since $g>1$, and that for any graph of groups $\G$ in normal
form associated to a handle-orbifold $\H$ one has $-\chi \ge 0$ unless
$\Gamma$ consists of a single vertex $v$, i.e. $\H$ consists just of a single 0-handle
orbifold $S^3/G_v$.

\medskip

Let $e$ be any edge of $\Gamma$ and denote by $a$ the order of its edge group; we will show
that  ${n \over a} \le 6(g-1).$

\medskip

Suppose first that $e$ is a closed edge (i.e., an edge that is a closed loop).  If $e$ is
the only edge of
$\G$ then
$$-\chi  \ge  {1 \over a} - {1 \over 2a} = {1 \over a}, \hskip 5mm
g-1 = -\chi n \ge {n \over 2a},  \hskip 5mm  {n \over a} \le 2(g-1).$$
If $e$ is closed and not the only edge then
$$-\chi  \ge  {1 \over a}, \hskip 5mm
g-1 = -\chi n \ge {n \over a},  \hskip 5mm  {n \over a} \le g-1.$$

Suppose that $e$ is not closed. If $e$ is the only edge of $\G$ then both
vertices of $e$ are isolated and
$$-\chi  \ge {1 \over a} - {1 \over 2a} - {1 \over 3a} = {1 \over 6a}, \hskip 5mm   g-1 =
-\chi \; n  \; \ge  \; {n \over 6a}, \hskip 5mm  {n \over a} \le 6(g-1).$$

\medskip

If $e$ is not closed, not the only edge and has exactly one isolated vertex then
$$-\chi  \ge {1 \over a} - {1 \over 2a} = {1 \over 2a}, \hskip 5mm   g-1 =
-\chi \; n  \; \ge  \; {n \over 2a}, \hskip 5mm  {n \over a} \le 2(g-1).$$
Finally, if $e$ is not closed, not the only edge and has no isolated vertex then
$$-\chi  \ge {1 \over a}, \hskip 5mm   g-1 =
-\chi \; n  \; \ge  \; {n \over a}, \hskip 5mm  {n \over a} \le g-1.$$

Concluding, in all cases we have
$${n \over a} \le 6(g-1).$$
In particular, if $\G$ has an edge whose edge group has order $a \le 60$ then
$n \le 6a(g-1) \le 360(g-1)$ which is linear in $g$; in particular, the quadratic bound of
Theorem 1 holds if $g \ge 15$.  So for the proof of the first part of the theorem we can
assume that all edge groups of $\G$ are cyclic or dihedral.

\medskip

Again, let $e$ be any edge of $\Gamma$; its edge group is either cyclic of order $a = b$ or
dihedral of order $a = 2b$. In each case this gives a cyclic subgroup $\Z_b$ of $G$ which
has a global fixed point in $H_g$. Let  $\K$ be a graph of groups in normal form associated
to the action of $\Z_b$ on $H_g$ (applying again Theorem 2), and $\chi' = \chi\K$.  Since
$\Z_b$ has a global fixed point in $H_g$, the graph of groups $\K$ must have a vertex $v$
with vertex group  $\Z_b$ such that, for the corresponding 0-handle orbifold $S^3/\Z_b$, the
action of the vertex group $\Z_b$ on the 0-handle $S^3$ has a global fixed point. Since no
cyclic subgroup of prime order can have two circles of fixed points, any nontrivial
subgroup of
$\Z_b$ has exactly the same circle of fixed points on the 0-handle $S^3$. Hence
all edges of $\K$ which have $v$ as a vertex have either trivial edge group or edge group
$\Z_b$ (by the Observation following Proposition 1).  We will show that $b \le 2g$.

\medskip

If all edges of $\K$ with vertex $v$ are closed then one must have trivial edge group
(since $\pi_1\K$ has no nontrivial finite normal subgroups), hence
$$-\chi'  \ge  1 - {1 \over b} =  {b-1 \over b}, \hskip 5mm
g-1 = -\chi' b \ge b-1,  \hskip 5mm  b \le g.$$

Otherwise there is a nonclosed edge $e'$ in $\K$ with vertex $v$, and the edge group of
$e'$ must be trivial (since $\K$ is in normal form, i.e. without trivial edges).

\medskip

If no vertex of $e'$ is isolated then
$$-\chi'  \ge  1, \hskip 5mm
g-1 = -\chi' b \ge  b,  \hskip 5mm  b \le g-1.$$
If exactly one vertex of $e'$ is isolated then
$$-\chi'  \ge  1 - {1 \over 2} =  {1 \over 2}, \hskip 5mm
g-1 = -\chi' b \ge {b \over 2},  \hskip 5mm  b \le 2(g-1).$$
If both vertices of $e'$ are isolated then
$$-\chi'  \ge 1  - {1 \over 2} - {1 \over b} = {b-2 \over 2b}, \hskip 5mm   g-1 =
-\chi' \; b  \; \ge  \; {b-2 \over 2}, \hskip 5mm  b \le 2g.$$

\medskip

In conclusion, in all cases we have
$$b \le 2g.$$

\medskip

Combining this with the inequality
$${n \over a} \le 6(g-1)$$
from above,  we obtain (since $a = b$ or $a = 2b$)
$$n = |G| = {n \over a} \cdot a  \le 6(g-1) \cdot 4g = 24 g(g-1).$$
This proves the first part of Theorem 1.

\bigskip

Now suppose that $G$ is a finite cyclic group of order $n$; then
all vertex and edge groups of $\G$ are  cyclic groups whose
orders divide $n$. Consider an edge $e$ of $\G$ with edge group $G_e \cong \Z_a$ whose order
$a$ realizes the minimum order over all edge groups.  If $e$ is the only edge of
$\G$ (with one or two distinct vertices) then $a=1$ (otherwise $\Z_a$ would be
a nontrivial finite normal subgroup of $\pi_1\G$), hence
$$-\chi  \ge  1- {1 \over 2} -
{1 \over 3} = {1 \over 6}, \hskip 1cm  g-1 = -\chi \; n  \ge  {n \over
6},  \hskip 1cm  n \le 6(g-1);$$
in particular, the quadratic bound of Theorem 1 in the cyclic case holds if $g \ge 3$.

\medskip

Suppose then that $e$ is not the only edge of $\Gamma$; consider another edge with
edge group $\Z_b$ which has a common vertex $v$ with $e$.  The vertex $v$
corresponds to a  0-handle $S^3/G_v$ of $\H$. The action of the
finite cyclic group $G_v$ on $S^3$ is orthogonal, so the union of the fixed
point sets of nontrivial elements of $G_v$ consists of at most two disjoint circles
$S^1$; also, a nontrivial subgroup of $G_v$ of prime order cannot fix two
circles (if we do not assume that the action of $G_v$ is orthogonal this follows from
Smith fixed point theory). Since by the Observation following Proposition 1 each
edge group is a maximal subgroup of a vertex group
$G_v$ fixing the point around which the corresponding 1-handle is attached, this
implies that either $a=b$ or $(a,b) = 1$ (the greatest common divisor).

\medskip

If $a=b$, since the edge group $G_e \cong
\Z_a$ is not normal in $\pi_1\G$ there must occur a situation of two edges
with a common vertex, one with edge group $\Z_a$ and the other with some edge group
$\Z_b$ such that $(a,b)=1$. This implies  $ab \le n$ and $a \le \sqrt {n}$ (since $a
\le b$).  If the considered edge with edge group $\Z_a$ has an isolated vertex (i.e., of
degree or valence 1), then
$$-\chi \;\; \ge \;\; {1 \over a} - {1 \over 2a} = {1 \over 2a} \;\; \ge \; \; {1 \over
2\sqrt {n}}, \hskip 7mm  g-1 = -\chi \; n \;\; \ge \;\; {n \over 2\sqrt {n}} = {\sqrt {n}
\over 2},$$
$$n \, \le \, 4(g-1)^2.$$
If the edge has no isolated vertex, one
obtains the stronger
$$-\chi \; \ge \; {1 \over a}  \; \ge \; {1 \over
\sqrt {n}}, \hskip 1cm g-1 \; \ge  \; \sqrt {n}, \hskip 1cm n \; \le \; (g-1)^2.$$

This proves Theorem 1 also in the cyclic case.

\medskip

Since the final assertion of Theorem 1 about the nonexistence of a linear bound
follows from the examples constructed in section 3, this completes the proof of
Theorem 1.

\bigskip \bigskip

{\bf Acknowledgment.}  The author was supported by a FRV grant from
Universit\`{a} degli Studi di Trieste.  He wants to thank the referee for his
various suggestions which helped to improve the paper.

\bigskip  \bigskip

\centerline {\bf References}

\medskip

\item {[B]} Z. Bao, {\it Maximum order of periodic outer automorphisms of a
free group,} J. Algebra 224  (2000), 437-453

\smallskip

\item {[C]}  M.J. Collins, {\it  On Jordan's theorem for complex linear groups,}
J. Group Theory 10   (2007),  411-423

\smallskip

\item {[Du]}  M.J. Dunwoody,  {\it  An equivariant sphere theorem,}
Bull. London Math. Soc. 17  (1985), 437-448

\smallskip

\item {[DV]}  P. Du Val, {\it  Homographies, Quarternions and Rotations,} Oxford
Math. Monographs, Oxford University Press 1964

\smallskip

\item {[HD]} A. Haefliger, Q.N. Du, {\it Appendice: une pr\'esentation du groupe
fondamental d'une orbifold,}  Ast\'erisque 115  (1984), 98-107

\smallskip

\item {[JR]}  W. Jaco, J.H. Rubinstein,  {\it  PL equivariant surgery and invariant
decompositions of 3-manifolds}  Adv. in Math. 73  (1989), 149-191

\smallskip

\item {[LN]} G. Levitt, J.-L. Nicolas, {\it On the maximum order of torsion
elements in $GL_n(\Bbb Z)$ and $Aut(F_n)$,} J. Algebra 208  (1998), 630-642

\smallskip

\item {[MMZ]} D. McCullough, A. Miller, B. Zimmermann,  {\it Group actions on
handlebodies,}  Proc. London Math. Soc.  59   (1989), 373-415

\smallskip

\item {[MSY]} W.H. Meeks, L. Simon, S.T. Yau,  {\it Embedded minimal surfaces, exotic
spheres, and manifolds with positive Ricci curvature,}  Ann. of Math. 116  (1982) , 621-659

\smallskip

\item {[MZ]} A. Miller, B. Zimmermann,  {\it  Large groups of symmetries of
handlebodies,}  Proc. Amer. Math. Soc. 106  (1989),  829-838

\smallskip

\item {[R]} F. Raymond,  {\it Classification of actions of the circle on 3-manifolds,}
Trans. Amer. Math. Soc. 131  (1968), 51-78

\smallskip

\item {[ScW]} P. Scott, T. Wall,  {\it  Topological methods in group theory,}
Homological Group Theory, London Math. Soc. Lecture Notes 36  (1979), Cambridge University
Press

\smallskip

\item {[Se]} J.P. Serre,  {\it  Trees,}  Springer, New York, 1980

\smallskip

\item {[WZ]} S. Wang, B. Zimmermann,  {\it  The maximum order finite groups of
outer automorphisms of free groups,}  Math. Z. 216  (1994), 83-87

\smallskip

\item {[Z1]} B. Zimmermann,  {\it \"Uber Abbildungsklassen von
Henkelk\"orpern,}  Arch. Math. 33  (1979),  379-382

\smallskip

\item {[Z2]} B. Zimmermann,  {\it \"Uber Hom\"oomorphismen n-dimensionaler Henkelk\"orper
und endliche Erweiterungen von Schottky-Gruppen,}  Comm. Math. Helv. 56   (1981), 474-486

\smallskip

\item {[Z3]} B. Zimmermann,  {\it Genus actions of finite groups on 3-manifolds,}
Michigan  Math. J.  43  (1996), 593-610

\item {[Z4]} B. Zimmermann,  {\it  Generators and relations for discontinuous groups,}
Generators and relations in Groups and Geometries,  NATO Advanced Study Institute Series
vol. 333  (1991), 407-436 (eds. Barlotti, Ellers, Plaumann, Strambach),  Kluwer Academic
Publishers

\smallskip

\item {[Z5]} B. Zimmermann,  {\it Finite groups of outer automorphism groups of free
groups,}  Glasgow Math. J. 38  (1996),  275-282

\bigskip \bigskip

Universit\`a degli Studi di Trieste

Dipartimento di Matematica e Geoscienze

34127 Trieste, Italy

\bye